\documentclass[10pt,journal,onecolumn]{IEEEtran}

\usepackage{bm, amssymb, amsmath, graphics, mathtools}
\newcommand{\T}{^{\mbox{\tiny T}}}
\newcommand{\B}[1]{{\bm #1}}
\newcommand{\ds}{\displaystyle}
\newcommand{\dd}{\; \text{d}}

\begin{document}

\title{The Theory of Connections. Connecting Points.}

\author{Daniele Mortari \\ \vspace{0.6cm} \emph{dedicated to John Lee Junkins}\footnote{Professor, Aerospace Engineering, Texas A\&M University, College Station, TX $77843$-$3141$, USA. E-mail: mortari$@$tamu.edu}}

\maketitle

\begin{abstract}
    This study introduces a procedure to obtain general expressions, $y = f(x)$, subject to linear constraints on the function and its derivatives defined at specified values. These \emph{constrained expressions} can be used describe functions with embedded specific constraints. The paper first shows how to express the most general explicit function passing through a single point in three distinct ways: linear, additive, and rational. Then, functions with constraints on single, two, or multiple points are introduced as well as those satisfying relative constraints. This capability allows to obtain general expressions to solve linear differential equations with no need to satisfy constraints (the ``subject to:'' conditions) as the constraints are already embedded in the \emph{constrained expression}. In particular, for expressions passing through a set of points, a generalization of the Waring's interpolation form, is introduced. The general form of additive \emph{constrained expressions} is introduced as well as a procedure to derive its coefficient functions, requiring the inversion of a matrix with dimensions as the number of constraints.
\end{abstract}

\section{Introduction}

This study introduces some special expressions, called \emph{constrained expressions}, to nonlinear functions subject to satisfy a set of linear constraints. These constraints are assigned at specified values of the independent variable ($x$) in term of the function and/or its derivatives. The general expression of any (nonlinear) function subject to pass through one point, two points, and multiple points, with single or multiple constraints, are provided. An application of these \emph{constrained expressions} is given in Ref. \cite{Mortari}, where it has been shown how to obtain least-squares solutions of linear nonhomogeneous differential equations of any order and nonconstant coefficients. This is done for initial and boundary values problems. Additional applications are in optimization as well as in path planning, just to name two. These \emph{constrained expressions} can be derived using \emph{any} set of functions spanning different function spaces. The resulting expression is provided in term of a new (nonlinear) function, $g (x)$, which is completely free to choose.

To show the most simple example of a \emph{constrained expression}, consider the following. The problem of writing the equation of \emph{all} linear functions passing through a specified point, $[x_1, y_1]$, is straightforward, $y (x) = m \, (x - x_1) + y_1$, with the line slope, $m$, free to choose. The focus of this paper is to write the equation of \emph{all} explicit functions passing through $[x_1, y_1]$. This includes \emph{all} functions: nonlinear, continuous, discontinuous, singular, periodic, etc.

Three different \emph{constrained expressions} are introduced, here. The first expression is a direct extension of the ``all-lines'' equation,
\begin{equation}\label{eq00}
    y (x) = p (x) \, (x - x_1) + y_1
\end{equation}
where $p (x)$ can be \emph{any} function satisfying, $p (x_1) \ne \infty$. The second expression, called \emph{additive}, is
\begin{equation}\label{eq01}
    y (x) = g (x) + [y_1 - g (x_1)] = g (x) + (y_1 - g_1),
\end{equation}
where $g (x_1) \ne \infty$, and the third expression, called \emph{rational}, is
\begin{equation}\label{eq02}
    y (x) = \dfrac{h (x)}{h (x_1)} \, y_1 = \dfrac{h (x)}{h_1} \, y_1,
\end{equation}
where $h (x)$ can be \emph{any} function satisfying $h (x_1) \ne 0$.

From Eq. (\ref{eq00}) and Eq. (\ref{eq01}) the following expression
\begin{equation}\label{eq00b}
    p (x) = \dfrac{g (x) - g_1}{x - x_1},
\end{equation}
is derived. This equation is interesting as $p (x)$ works as the derivative of $g (x)$, but for finite variations.

The expressions introduced in Eqs. (\ref{eq00},\ref{eq01},\ref{eq02}) can be used to represent \emph{all} explicit functions passing through the point $[x_1, y_1]$, with the only exception to those satisfying $p (x_1) = g (x_1) = \infty$ and $h (x_1) = 0$. It is possible also to combine Eqs. (\ref{eq00},\ref{eq01},\ref{eq02}) to obtain
\begin{eqnarray*}
    y (x) &=& p (x) \, (x - x_1) + \dfrac{h (x)}{h (x_1)} \, y_1 \\
    y (x) &=& p (x) \, (x - x_1) + g (x) + (y_1 - g_1) \\
    y (x) &=& g (x) + \dfrac{h (x)}{h (x_1)} \, (y_1 - g_1) \\
    y (x) &=& p (x) \, (x - x_1) + g (x) + \dfrac{h (x)}{h (x_1)} \, (y_1 - g_1)
\end{eqnarray*}

This study investigates how to derive these expressions subject to a variety of linear constraints, such as functions passing through multiple points with assigned derivatives, or subject to multiple relative constraints, as well as periodic functions subject to multiple point constraints. Potential applications for these types of expressions are many. For example, they can be used in optimization problems using expressions with the assigned constraints already embedded. This is particularly useful when solving linear differential equations, especially for boundary values problems, where the integration techniques become complicate to satisfy the constraints. In the contrary, using the expressions provided in this study, linear differential equations can be rewritten for functions with already embedded the constraints (a.k.a., the ``subject to'' conditions). This approach provides to the initial and boundary values problems a unified solving framework. Specifically, for linear differential equations, a least-squares approach solving technique \cite{Mortari} has been developed with order of magnitudes accuracy gain with respect the step-varying Runge-Kutta-Fehlberg method.

Particularly important is the fact that Eq. (\ref{eq01}), and many other equations provided in this study, can be immediately extended to vectors (bold lower cases), to matrices (upper cases), and to higher dimensional tensors. Specifically, for vectors and matrices, Eq. (\ref{eq01}) becomes
\begin{equation*}
    \B{y} (x) = \B{g} (x) + (\B{y}_1 - \B{g}_1) \qquad \text{and} \qquad Y (x) = G (x) + (Y_1 - G_1).
\end{equation*}

In the following sections this study derives \emph{constrained expressions} satisfying:
\begin{itemize}
  \item Constraints in one point;
  \item Constraints in two points and extension to $n$ points;
  \item $m$ constraints in $n \le m$ points;
  \item Single and multiple relative constraints;
  \item Constraints on continuous and discontinuous periodic functions.
\end{itemize}

\section{Constraints on one point}

Functions subject to, $y (x_1) = y_1$, can be derived using the general form
\begin{equation}\label{eq03}
    y (x) = g (x) + \xi_1 \, h (x)
\end{equation}
where $\xi_1$ is an unknown coefficient and $g (x)$ and $h (x)$ can be \emph{any} two functions, satisfying $g (x_1) = g_1 \ne \infty$ and $h (x_1) = h_1 \ne 0$. The constraint, $y (x_1) = y_1$, allows us to derive the expression of $\xi_1$, getting,
\begin{equation}\label{eq04}
     y (x) = g (x) + \dfrac{h (x)}{h_1} \, (y_1 - g_1),
\end{equation}
which is a combination of Eq. (\ref{eq01}) and Eq. (\ref{eq02}). In particular, if $g (x) = 0$ or if $h (x) = 1$ we obtain the expressions already provided in Eq. (\ref{eq01}) and Eq. (\ref{eq02}), respectively.

\subsection{Constraint on one derivative}

Equation (\ref{eq03}) can also be used to obtain a function satisfying the constraint on the $n$-th derivative, $\left.\dfrac{\dd^n y}{\dd x^n}\right|_{x = x_1} = y_1^{(n)}$, getting
\begin{equation*}
     y (x) = g (x) + \dfrac{h (x)}{h_1^{(n)}} \, \left[y_1^{(n)} - g_1^{(n)}\right]
\end{equation*}
where $g_1^{(n)} \ne \infty$ and $h_1^{(n)} \ne 0$. In particular, if $h (x) = x^n$, then $y (x) = g (x) + \dfrac{x^n}{n!} \, \left[y_1^{(n)} - g_1^{(n)}\right]$.

\subsection{Constraints on two derivatives}

Functions constrained by $\left.\dfrac{\dd^p y}{\dd x^p}\right|_{x = x_1} = y_1^{(p)}$ and $\left.\dfrac{\dd^q y}{\dd x^q}\right|_{x = x_1} = y_1^{(q)}$, where $0 \le p < q$, can be derived using the expression
\begin{equation}
    y (x) = g (x) + \xi_p \, h_p (x) + \xi_q \, h_q (x).
\end{equation}
The two constraints allow to compute the coefficients $\xi_p$ and $\xi_p$, by
\begin{equation*}
    \begin{Bmatrix} y^{(p)}_1 - g^{(p)}_1 \\ y^{(q)}_1 - g^{(q)}_1\end{Bmatrix} = \begin{bmatrix} h^{(p)}_{p1} & h^{(p)}_{q1}\\ h^{(q)}_{p1} & h^{(q)}_{q1}\end{bmatrix} \begin{Bmatrix} \xi_p\\ \xi_q\end{Bmatrix}
\end{equation*}
that is
\begin{equation*}
    \begin{Bmatrix} \xi_p\\ \xi_q\end{Bmatrix} = \dfrac{1}{h^{(p)}_{p1} h^{(q)}_{q1} - h^{(p)}_{q1} h^{(q)}_{p1}} \begin{bmatrix} h^{(q)}_{q1} & -h^{(p)}_{q1}\\ -h^{(q)}_{p1} & h^{(p)}_{p1}\end{bmatrix} \begin{Bmatrix} y^{(p)}_1 - g^{(p)}_1\\ y^{(q)}_1 - g^{(q)}_1\end{Bmatrix}.
\end{equation*}
Solution exists as long as $h^{(p)}_{p1} h^{(q)}_{q1} \ne h^{(p)}_{q1} h^{(q)}_{p1}$. Under this assumption the searched constrained expression is
\begin{equation*}
    y (x) = g (x) + \dfrac{h^{(q)}_{q1} h_p (x) - h^{(q)}_{p1} h_q (x)}{h^{(p)}_{p1} h^{(q)}_{q1} - h^{(p)}_{q1} h^{(q)}_{p1}} \left[y^{(p)}_1 - g^{(p)}_1\right] + \dfrac{h^{(p)}_{p1} h_q (x) - h^{(p)}_{q1} h_p (x)}{h^{(p)}_{p1} h^{(q)}_{q1} - h^{(p)}_{q1} h^{(q)}_{p1}} \left[y^{(q)}_1 - g^{(q)}_1\right]
\end{equation*}
The relationship $h^{(p)}_{p1} h^{(q)}_{q1} \ne h^{(p)}_{q1} h^{(q)}_{p1}$ provides us the following considerations. In order to admit solution, $h_p (x)$ and $h_q (x)$ must have a nonzero $p$ and $q$ derivatives, $h^{(p)}_p (x) \ne 0$ and $h^{(q)}_q (x) \ne 0$, respectively. For example, if $p = 0$ and $q = 3$ then $h_p (x)$ must be (at least) a nonzero constant and $h_q (x)$ must be (at least) cubic, $h_q (x) = x^3$. For instance, this requirement can always be obtained by setting $h_p (x) = \dfrac{x^p}{p!}$ and $h_q (x) = \dfrac{x^q}{q!}$. Then, $\xi_p$ and $\xi_q$ are derived from the constraints,
\begin{equation*}
    \begin{Bmatrix} y^{(p)}_1 - g^{(p)}_1\\ y^{(q)}_1 - g^{(q)}_1\end{Bmatrix} = \begin{bmatrix} 1 & \dfrac{x_1^{(q - p)}}{(q - p)!}\\ 0 & 1\end{bmatrix} \begin{Bmatrix} \xi_p \\ \xi_q\end{Bmatrix} \quad \rightarrow \quad \begin{Bmatrix} \xi_p \\ \xi_q\end{Bmatrix} = \begin{bmatrix} 1 & -\dfrac{x_1^{(q - p)}}{(q - p)!} \\ 0 & 1\end{bmatrix} \begin{Bmatrix} y^{(p)}_1 - g^{(p)}_1\\ y^{(q)}_1 - g^{(q)}_1\end{Bmatrix}
\end{equation*}
and the solution is
\begin{equation}\label{eq05}
    y (x) = g (x) + \dfrac{x^p}{p!} \left[y^{(p)}_1 - g^{(p)}_1\right] + \left[\dfrac{x^{(q - p)}}{q!} - \dfrac{x_1^{(q - p)}}{p!(q - p)!} \right] \, x^p \left[y^{(q)}_1 - g^{(q)}_1\right].
\end{equation}
In the common case where function ($p = 0$) and first derivative ($q = 1$) are assigned in one point, the Eq. (\ref{eq05}) becomes
\begin{equation}\label{eqmin}
    y (x) = g (x) + (y_1 - g_1) + (x - x_1) (\dot{y}_1 - \dot{g}_1).
\end{equation}
Note that, if $g (x) = a$ or $g(x) = a + b \, x$, then the previous equation reduces to $y (x) = y_1 + \dot{y}_1 (x - x_1)$. The reason is, to obtain Eq. (\ref{eqmin}), we have implicitly selected $h_p (x) = 1$ and $h_q (x) = x$. Therefore, in order for $g (x)$ to provides additional variation, then its expression must be \emph{at least quadratic}, that is, at least one degree greater than $h_p (x)$ and $h_q (x)$.

\subsection{Constraints on $n$ derivatives}

The constrained expression when the function and its first $n$ derivatives are assigned in one point can be written as
\begin{equation}\label{Taylor}
    \boxed{ y (x) = g (x) + \ds\sum_{k = 0}^n \dfrac{(x - x_1)^k}{k!} \left(\dfrac{\dd^k y}{\dd x^k} - \dfrac{\dd^k g}{\dd x^k}\right)_{x_1} }
\end{equation}
Equation (\ref{Taylor}) satisfies all the constraints, $y (x_1) = y_1$, $\dot{y} (x_1) = \dot{y}_1$, $\ddot{y} (x_1) = \ddot{y}_1$, $\dddot{y} (x_1) = \dddot{y}_1$, and so on. In particular, for $n = \infty$, this equation becomes the combination of two Taylor series, expanded around $x_1$, one for $y (x)$ and the other for $g (x)$, respectively,
\begin{equation}\label{Taylor2}
    y (x) = \ds\sum_{k = 0}^{\infty} \left.\dfrac{\dd^k y}{\dd x^k}\right|_{x_1} \dfrac{(x - x_1)^k}{k!} \qquad \text{and} \qquad g (x) = \ds\sum_{k = 0}^{\infty} \left.\dfrac{\dd^k g}{\dd x^k}\right|_{x_1} \dfrac{(x - x_1)^k}{k!}.
\end{equation}
Therefore, as $n \to \infty$ the set of potential functions that can be described by Eq. (\ref{Taylor}) converges to a single function only, the $y (x)$ defined in Eq. (\ref{Taylor2}). This means that, as $n$ increases the variations affecting $y (x)$, due to the $g (x)$ variations, decrease. Note that, Eq. (\ref{Taylor}) is linear in $g (x)$ and its derivatives and, therefore, it may be possible to take advantage of $g (x)$ and increase the predictions based on truncated Taylor serie.

For vectors made of $m$ \emph{independent variables} such as $\B{y}\T (x) = \{y_1 (x), \; y_2 (x), \; \cdots, \; y_m (x) \}$, and subject to $n$ constraints at $x = x_1$, Eq. (\ref{Taylor}) become
\begin{equation}\label{vector4}
    \B{y} (x) = \B{g} (x) + \ds\sum_{k = 0}^n \dfrac{(x - x_1)^k}{k!} \left(\dfrac{\dd^k \B{y}}{\dd x^k} - \dfrac{\dd^k \B{g}}{\dd x^k}\right)_{x_1}
\end{equation}
In this case the $\B{g} (x)$ vector can be expressed as linear combination of $n$ common basis functions, $\B{h}\T (x) = \{h_1 (x), \; h_2 (x), \; \cdots, \; h_n (x)\}$,
\begin{equation}\label{vector5}
    \B{g} (x) = \begin{bmatrix} \B{\xi}_1, & \B{\xi}_2, & \cdots, & \B{\xi}_n\end{bmatrix}\T \B{h} (x) = \Xi \, \B{h} (x) \qquad \rightarrow \qquad \left.\dfrac{\dd^k \B{g}}{\dd x^k}\right|_{x_1} = \Xi \, \left.\dfrac{\dd^k \B{h}}{\dd x^k}\right|_{x_1}
\end{equation}
where $\Xi$ is a $n\times m$ matrix of the unknown coefficients. Then, Eq. (\ref{vector4}) becomes
\begin{equation}\label{vector6}
    \B{y} (x) = \Xi \left[\B{h} (x) - \ds\sum_{k = 0}^n \dfrac{(x - x_1)^k}{k!} \B{h}^{(k)}_1 \right] + \ds\sum_{k = 0}^n \dfrac{(x - x_1)^k}{k!} \B{y}^{(k)}_1
\end{equation}

Equation (\ref{Taylor}) has interesting applications for vectors made of subsequent derivatives, such as, $\B{y}_d\T (x) = \{y, \, \dot{y}, \, \ddot{y}, \, \dddot{y}, \cdots\}$, where the subscript ``$d$'' identifies this specific kind of vector. The constrained expression for this vector (which often appear in dynamics) is,
\begin{equation}\label{vector1}
    \B{y}_d (x) = \begin{Bmatrix} y (x) \\ \dot{y} (x) \\ \ddot{y} (x) \\ \vdots\end{Bmatrix} = \B{g}_d + B (x, x_0) (\B{y}_{d0} - \B{g}_{d0}) = \begin{Bmatrix} g (x) \\ \dot{g} (x) \\ \ddot{g} (x) \\ \vdots\end{Bmatrix} + B (x, x_0) \begin{Bmatrix} y_0 - g_0 \\ \dot{y}_0 - \dot{g}_0 \\ \ddot{y}_0 - \ddot{g}_0 \\ \vdots\end{Bmatrix}
\end{equation}
where the $B (x, x_0)$ matrix is
\begin{equation}\label{vector2}
    B (x, x_0) = \begin{bmatrix} 1 & (x - x_0) & \dfrac{(x - x_0)^2}{2!} & \dfrac{(x - x_0)^3}{3!} & \dfrac{(x - x_0)^4}{4!} & \cdots \\ 0 & 1 & (x - x_0) & \dfrac{(x - x_0)^2}{2!} & \dfrac{(x - x_0)^3}{3!} & \cdots \\ 0 & 0 & 1 & (x - x_0) & \dfrac{(x - x_0)^2}{2!} & \cdots \\ \vdots & \vdots & \vdots & \vdots & \vdots & \ddots\end{bmatrix}
\end{equation}
whose elements ($i$-row and $j$-column) can be defined as
\begin{equation*}
    B (i, j) = \left\{\begin{array}{lll} = 0 & \text{if} & i > j \\ = 1 & \text{if} & i = j \\ = \dfrac{(x - x_0)^{j-i}}{(j-i)!} & \text{if} & i < j\end{array}\right. \quad \text{and} \quad \dot{B} (i, j) = \left\{\begin{array}{lll} = 0 & \text{if} & i \ge j \\ = 1 & \text{if} & i = j - 1 \\ = \dfrac{(x - x_0)^{j-i-1}}{(j-i-1)!} & \text{if} & i < j - 1\end{array}\right.
\end{equation*}
To use this kind of constrained equations, the $\B{g}_d (x)$ vector can be expressed as a linear combination of a set of $n$ basis functions, $\B{h}\T (x) = \{ h_1 (x), \; h_2 (x), \; \cdots, \; h_n (x) \}$,
\begin{equation}\label{vector3}
    \B{g}_d = \begin{bmatrix} \B{h} (x), & \dot{\B{h}} (x), & \ddot{\B{h}} (x), & \dddot{\B{h}} (x), & \cdots\end{bmatrix}\T \, \B{\xi} = H \, \B{\xi}
\end{equation}
where $\B{\xi}$ is a vector of $n$ unknown coefficients to be found.  Note that, when the vector is made of independent variables, the number of unknown coefficients is higher: all elements of the $\Xi$ matrix.

\section{Constraints in $n$ points}

Waring polynomials \cite{Waring}, better known as ``Lagrange polynomials,'' are used for polynomial interpolation. This interpolation formula was first discovered in 1779 by Edward Waring, then rediscovered by Euler in 1783, and then published by Lagrange in 1795 \cite{Michiel}. A two-points Waring polynomial is the equation of a line passing through two points, $[x_1, y_1]$ and $[x_2, y_2]$,
\begin{equation*}
    y (x) =  y_1 \, \dfrac{x - x_2}{x_1 - x_2} + y_2 \, \dfrac{x - x_1}{x_2 - x_1}
\end{equation*}
while the general $n$-points Waring polynomial is the ($n-1$) degree curve passing through $n$ points,
\begin{equation*}
    y (x) = \ds\sum_{k = 1}^n y_k \ds\prod_{i \ne k} \dfrac{x - x_i}{x_k - x_i}.
\end{equation*}

Inspired by Waring polynomials the expression of a function subject to pass through two or more distinct points can be provided. In fact, Waring polynomials can be generalized to obtain \emph{all} nonlinear functions connecting $[x_1, y_1]$ and $[x_2, y_2]$ using additive expression,
\begin{equation}\label{eqw1}
    y (x) = g(x) + \dfrac{x - x_2}{x_1 - x_2} \, (y_1 - g_1) + \dfrac{x - x_1}{x_2 - x_1} \, (y_2 - g_2),
\end{equation}
where $g_1\ne \infty$ and $g_2\ne \infty$. Equation (\ref{eqw1}) describes all nonlinear functions passing through two points. The generalization is immediate. The expression representing \emph{all} functions passing through a set of $n$ points can be generalized as
\begin{equation}\label{eqw4}
    \boxed{ y (x) = g(x) + \ds\sum_{k = 1}^n (y_k - g_k) \ds\prod_{i\ne k} \dfrac{x - x_i}{x_k - x_i}. }
\end{equation}
Equation (\ref{eqw4}) extends the Waring's interpolation form using $n$ points to any function satisfying $g (x_k) \ne \infty$, where $k = 1, \cdots, n$.

\begin{figure}[ht]
    \centering\includegraphics[width=\linewidth]{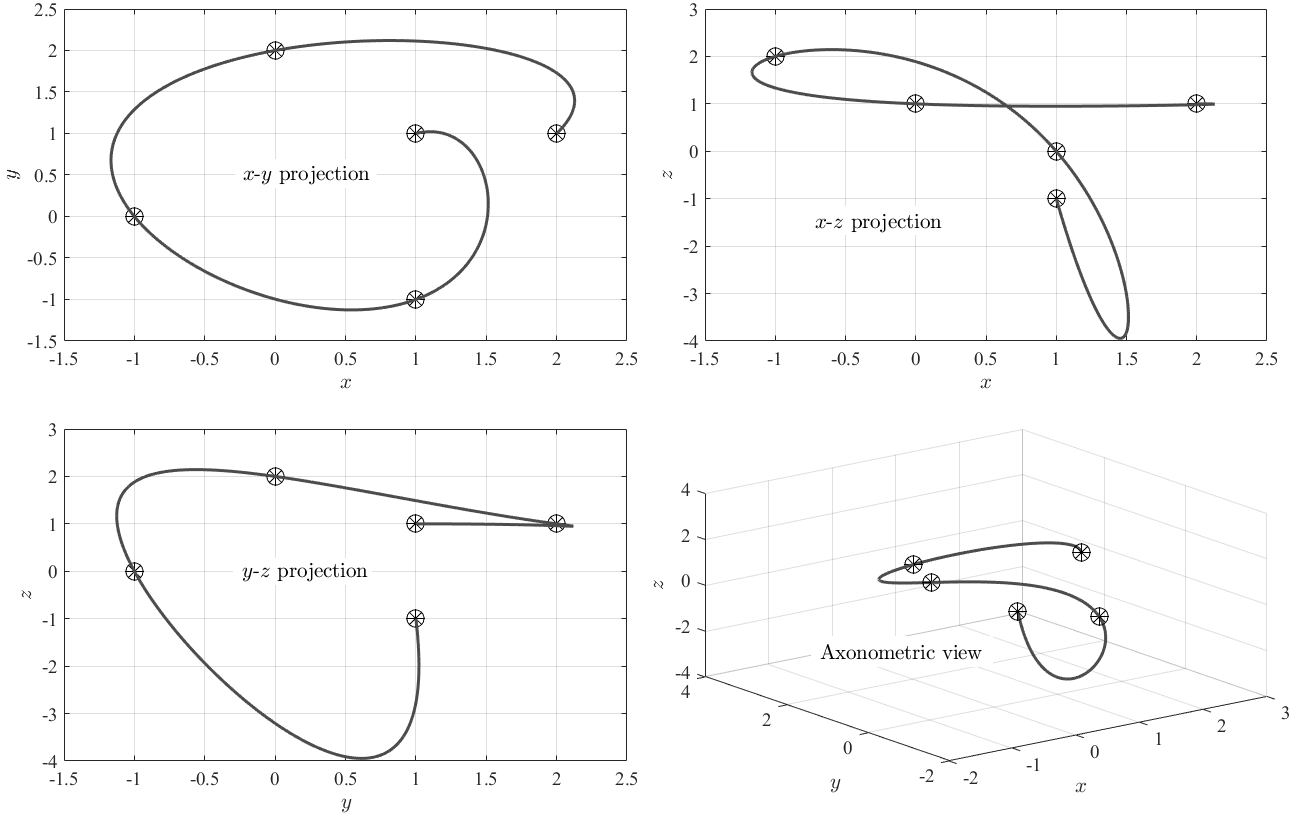}
	\caption{Trajectory obtained using $t\in[0,1]$ and data given in Table \ref{t1}}
    \label{vector}
\end{figure}

In particular, for a time-varying $n$-dimensional vector, $\B{y}$, passing through a set of $m$ points, $[\B{y}_1, \, \B{y}_2, \, \cdots, \, \B{y}_m]$, at respective times, $[t_1, \, t_2, \, \cdots, \, t_m]$, Eq. (\ref{eqw4}) becomes
\begin{equation*}
    \B{y} (t) = \B{g} (t) + \ds\sum_{k = 1}^m (\B{y}_k - \B{g}_k) \ds\prod_{i\ne k} \dfrac{t - t_i}{t_k - t_i}.
\end{equation*}
Using the five points given in Table \ref{t1} and $\B{g} = \{\, \sin t, \; e^t, \; 1 - t^2 \,\}\T$
the trajectory shown in Fig. \ref{vector} is obtained

\begin{table}[ht]
  \centering
  \label{t1}
  \begin{tabular}{c||ccrrr}
    ~ & $\B{y}_1$ & $\B{y}_2$ & $\B{y}_3$ & $\B{y}_4$ & $\B{y}_5$ \\ \hline\hline
    $x$ & 2 & 0 & $-1$ & 1 & 1 \\
    $y$ & 1 & 2 & 0 & $-1$ & 1 \\
    $z$ & 2 & 1 & 2 & 0 & $-1$ \\
    \hline\hline
  \end{tabular}
\end{table}
\vspace{-0.5cm}
\begin{center} Table 1. Selected points coordinates.\end{center}

Next two subsections contain two examples of two-point function constraints. In the first example function and first derivative are assigned in two distinct points while in the second example the first derivatives are assigned in two distinct points.

\subsection{Two constraints example \#1}

Consider a function subject to: $y (x_1) = y_1$ and $\dot{y} (x_2) = \dot{y}_2$. This function can be expressed as
\begin{equation*}
    y (x) = g (x) + \xi_0 \, h_0 (x) + \xi_1 \, h_1 (x).
\end{equation*}
In particular, if $h_0 (x) = 1$ and $h_1 (x) = x$, then $y (x) = g (x) + \xi_0 + \xi_1 \, x$, where $\xi_0$ and $\xi_1$ are two constants and $g (x)$ can be any nonlinear function. The two constraints imply solving the system
\begin{equation*}
    \begin{Bmatrix} y_1 - g_1\\ \dot{y}_2 - \dot{g}_2\end{Bmatrix} = \begin{bmatrix} 1 & x_1 \\ 0 & 1\end{bmatrix} \begin{Bmatrix} \xi_0 \\ \xi_1\end{Bmatrix}
\end{equation*}
from which we derive the solution
\begin{equation}\label{c1}
    y (x) = g (x) + (y_1 - g_1) + (x - x_1) \, (\dot{y}_2 - \dot{g}_2).
\end{equation}
Again, note that, if $g(x) = a + b \, x$, then Eq. (\ref{c1}) reduces to $y (x) = y_1 + \dot{y}_2 (x - x_1)$, that is, an expression where $g (x)$ disappears, unaffecting the constrained expression.

\subsection{Two constraints example \#2}

Consider a function subject to: $\dot{y} (x_1) = \dot{y}_1$ and $\dot{y} (x_2) = \dot{y}_2$. This function can be expressed as,
\begin{equation*}
    y (x) = g (x) + \xi_0 \, h_0 (x) + \xi_1 \, h_1 (x).
\end{equation*}
Assuming, $h_0 (x) = x$ and $h_1 (x) = x^2$, then, $y (x) = g (x) + \xi_0 \, x + \xi_1 \, x^2$, where $\xi_0$ and $\xi_1$ are two constants and $g (x)$ can be any nonlinear function. In this example we see another property affecting the $h_k (x)$ functions. In addition to the property that they \emph{cannot span the same function space}, all functions $h_k (x)$ \emph{must admit a nonzero derivative with same order of the constraint's minimum derivative order}. This means we cannot adopt for these constraints, as in the previous example, $h_0 (x) = 1$ and $h_1 (x) = x$.

The two constraints imply solving the system
\begin{equation*}
    \begin{Bmatrix} \dot{y}_1 - \dot{g}_1\\ \dot{y}_2 - \dot{g}_2\end{Bmatrix} = \begin{bmatrix} 1 & 2 x_1 \\ 1 & 2 x_2\end{bmatrix} \begin{Bmatrix} \xi_0 \\ \xi_1\end{Bmatrix}
\end{equation*}
giving the solution,
\begin{equation}\label{c2}
    y (x) = g (x) + \dfrac{x (2 x_2 - x)}{2 (x_2 - x_1)}  (\dot{y}_1 - \dot{g}_1) + \dfrac{x (x - 2 x_1)}{2 (x_2 - x_1)} (\dot{y}_2 - \dot{g}_2).
\end{equation}
Note that, if $g (x) = b \, x + c \, x^2$, then Eq. (\ref{c2}) reduces to the same equation, independently what $b$ and $c$ are. This means that, in order for $g (x)$ to play a role in the constrained expression, $g(x) \ne c_0 \, h_0 (x) + c_1 \, h_1 (x)$.

\subsection{Four constraints example}

Consider to find a \emph{constrained expression} subject to four constraints. Assume these constraints be specified in the following 3 distinct points,
\begin{equation}\label{eq07}
    \left.\dfrac{\dd^2 y}{\dd x^2}\right|_{\B{x}_1 =-1} = \ddot{y}_{\B{x}_1}, \quad y (\B{x}_2 = 0) = y_{\B{x}_2}, \quad y (\B{x}_3 = 2) = y_{\B{x}_3}, \quad \text{and} \quad \left.\dfrac{\dd y}{\dd x}\right|_{\B{x}_4 = 2} = \dot{y}_{\B{x}_4}.
\end{equation}
Let us find the constrained expression using monomials
\begin{equation}\label{eq08}
    y (x) = g (x) + \xi_0 + \xi_1 \, x + \xi_2 \, x^2 + \xi_3 \, x^3
\end{equation}
The constraints imply
\begin{equation*}
    \begin{Bmatrix} \ddot{y}_{\B{x}_1} - \ddot{g}_{\B{x}_1}\\ y_{\B{x}_2} - g_{\B{x}_2}\\ y_{\B{x}_3} - g_{\B{x}_3}\\ \dot{y}_{\B{x}_4} - \dot{g}_{\B{x}_4}\end{Bmatrix} = \begin{bmatrix} 0 & 0 & 2 & 6 \, \B{x}_1 \\ 1 & \B{x}_2 & \B{x}_2^2 & \B{x}_2^3\\ 1 & \B{x}_3 & \B{x}_3^2 & \B{x}_3^3\\ 0 & 1 & 2 \, \B{x}_4 & 3 \, \B{x}_4^2\end{bmatrix} \begin{Bmatrix} \xi_0\\ \xi_1\\ \xi_2\\ \xi_3\end{Bmatrix} = \begin{bmatrix*}[r] 0 & 0 & 2 & -6 \\ 1 & 0 & 0 & 0\\ 1 & 2 & 4 & 8\\ 0 & 1 & 4 & 12\end{bmatrix*} \begin{Bmatrix} \xi_0\\ \xi_1\\ \xi_2\\ \xi_3\end{Bmatrix}
\end{equation*}
providing the expressions for the four coefficients $\xi_k$
\begin{equation*}
    \begin{Bmatrix} \xi_0\\ \xi_1\\ \xi_2\\ \xi_3\end{Bmatrix} = \dfrac{1}{28} \begin{bmatrix*}[r] 0 & 28 & 0 & 0\\ -8 & -24 & 24 & -20\\ 8 & 3 & -3 & 6\\ -2 & 1 & -1 & 2\end{bmatrix*} \begin{Bmatrix} \ddot{y}_{\B{x}_1} - \ddot{g}_{\B{x}_1}\\ y_{\B{x}_2} - g_{\B{x}_2}\\ y_{\B{x}_3} - g_{\B{x}_3}\\ \dot{y}_{\B{x}_4} - \dot{g}_{\B{x}_4}\end{Bmatrix}
\end{equation*}
Then, substituting these expressions in Eq. (\ref{eq08}) we obtain the searched constrained expression,
\begin{equation}\label{c3}
\begin{split}
    y (x) = g (x) & + \dfrac{- 8 \, x + 8 \, x^2 - 2 \, x^3}{28} \, (\ddot{y}_{\B{x}_1} - \ddot{g}_{\B{x}_1}) + \dfrac{28 - 24 \, x +  3 \, x^2 + x^3}{28} \, (y_{\B{x}_2} - g_{\B{x}_2}) + \\ ~ & + \dfrac{24 \, x -  3 \, x^2 - x^3}{28} \, (y_{\B{x}_3} - g_{\B{x}_3}) + \dfrac{- 20 \, x +  6 \, x^2 + 2 \, x^3}{28} \, (\dot{y}_{\B{x}_4} - \dot{g}_{\B{x}_4}),
\end{split}
\end{equation}
satisfying all constraints defined in Eq. (\ref{eq06}). Again, the free function must satisfy $g (x) \ne c_0 + c_1 \, x + c_2 \, x^2 + c_3 \, x^3$, otherwise the minimal constrained equation,
\begin{equation*}
    y (x) = \dfrac{- 8 \, x + 8 \, x^2 - 2 \, x^3}{28} \, \ddot{y}_{\B{x}_1} + \dfrac{28 - 24 \, x +  3 \, x^2 + x^3}{28} \, y_{\B{x}_2} + \dfrac{24 \, x -  3 \, x^2 - x^3}{28} \, y_{\B{x}_3} + \dfrac{- 20 \, x +  6 \, x^2 + 2 \, x^3}{28} \, \dot{y}_{\B{x}_4},
\end{equation*}
is obtained.

\subsection{Issues using monomials}

The solution to the general problem of $n$ constraints on $m$ points exists as long as the matrix to derive the $\xi_k$ coefficients is not singular. For instance, using the following constraints,
\begin{equation*}
    \left.\dfrac{\dd^3 y}{\dd x^3}\right|_{\B{x}_1} = y^{(3)}_{\B{x}_1}, \qquad y (\B{x}_2) = y_{\B{x}_2}, \qquad \left.\dfrac{\dd^3 y}{\dd x^3}\right|_{\B{x}_2} = y^{(3)}_{\B{x}_2}, \qquad \text{and} \qquad \left.\dfrac{\dd^3 y}{\dd x^3}\right|_{\B{x}_3} = y^{(3)}_{\B{x}_3},
\end{equation*}
the monomial formalism of Eq. (\ref{eq08}) cannot be adopted because the corresponding coefficient matrix,
\begin{equation*}
    \begin{bmatrix} 0 & 0 & 0 & 6 \\ 1 & \B{x}_2 & \B{x}_2^2 & \B{x}_2^3 \\ 0 & 0 & 0 & 6 \\ 0 & 0 & 0 & 6\end{bmatrix},
\end{equation*}
has rank 2 and cannot be inverted. In this case, the minimal monomial that can be used in Eq. (\ref{eq08}) must be $y (x) = g (x) + \xi_0 + \xi_1 \, x^3 + \xi_2 \, x^4 + \xi_3 \, x^5$. It is easy to prove by induction that any smooth function $f$ for which $\dd^n f/\dd x^n = 0$ is a polynomial of degree less than $n$ and the only smooth functions for which a derivative of some order is identically zero are polynomials.

The use of monomials (which usually led to simple constrained expressions) can still be adopted for simple constrained expressions, while for the most general case) and to avoid singular coefficient matrix), the selection of a set of functions that are infinitely differentiable, such as exponentials, logarithm, trigonometric functions, rational functions, etc., are preferred.

In the next section a general method to derive \emph{constrained expressions} is provided for the general case of functions with $n$ constraints on $m$ points.

\section{Coefficients functions of constrained expressions}

Consider the general case of a function with $n$ distinct constraints, $\left.\dfrac{\dd^{\B{d}_k} y}{\dd x^{\B{d}_k}} \right|_{\B{x}_k} = y_{\B{x}_k}^{(\B{d}_k)}$, where the $n$-element vector, $\B{d}$, contains the derivative orders (of the $n$ constraints) and the $n$-element vector, $\B{x}$, indicates where the constraints are specified. For instance, for the constraints specified in Eq. (\ref{eq06}), $\B{d}$ and $\B{x}$ vectors are $\B{d} = \{2, \; 0, \; 0, \; 1\}\T$ and $\B{x} = \{-1, \; 0, \; 2, \; 2\}\T$, respectively.

In this section, we derive a different approach to derive constrained expressions. This is motivated by the fact that Eqs. (\ref{eq05}), (\ref{eqw1}), (\ref{c1}), (\ref{c2}), and (\ref{c3}) all share the same formalism, with a number of terms equal to the number of constraints, expressed as $\left(\dfrac{\dd^{\B{d}_k} y}{\dd x^{\B{d}_k}} - \dfrac{\dd^{\B{d}_k} g}{\dd x^{\B{d}_k}}\right)_{\B{x}_k}$, and multiplied by some coefficient functions assuming values of 1 if computed at the constraint's coordinate value while all the remaining coefficient functions are 0, and viceversa. This property suggests to search the constrained expressions as,
\begin{equation}\label{eq06}
    \boxed{ y (x) = g (x) + \ds\sum_{k = 1}^n \beta_k (x, \B{x}) \left[y_{\B{x}_k}^{(\B{d}_k)} - g_{\B{x}_k}^{(\B{d}_k)}\right] \quad \text{where} \quad \beta_i^{(\B{d}_k)} (x_k, \B{x}) = \delta_{ki}, }
\end{equation}
where $\delta_{ki}$ is the Kronecker and $n$ the number of constraints. The \emph{coefficient functions}, $\beta_k (x, \B{x})$, of this expression are such that when the $k$-th constraint, $\left.\dfrac{\dd^{\B{d}_k} y}{\dd x^{\B{d}_k}}\right|_{\B{x}_k} = y_{\B{x}_k}^{(\B{d}_k)}$, is verified, then $\beta_k^{(\B{d}_k)} (x_k, \B{x}) = 1$, while all the other coefficient functions, $\beta_i^{(\B{d}_k)} (x_k, \B{x}) = 0$, are zeros ($i\ne k$). Given a set of constraints, this property allows us to find the $\beta_k (x, \B{x})$ expressions. This is explained in the next section.

\subsection{Coefficients functions derivation}

As the number of constraints ($m$) increases, the approach previously proposed to find constrained expressions becomes complicate with the risk of obtaining singular matrix when computing the coefficient vector, $\B{\xi}$. For this reason this section propose a new procedure to compute all the $\beta_k$ functions at once, provided that the $\beta_k$ functions are expressed as scalar product of a set of $m$ linearly \emph{independent} functions (preferably indefinitely differentiable), contained in the vector $\B{h} (x) = \{h_1 (x), h_2 (x), \cdots, h_m (x)\}\T$ (where $m \le n$), and the vectors of the unknown functions coefficients,
\begin{equation*}
    \beta_k (x) = \B{h}\T (x) \, \B{\xi}_k.
\end{equation*}
To be clear, the procedure is shown using the following $n = 4$ constraints example,
\begin{equation*}
    y (x_1) = y_1, \qquad \dddot{y} (x_1) = \dddot{y}_1, \qquad \dot{y} (x_2) = \dot{y}_2, \qquad \text{and} \qquad y (x_3) = y_3.
\end{equation*}
To compute the function $\beta_1 (x)$ associated to the first constraint, $y (x_1) = y_1$, the following relationships,
\begin{equation*}
    \beta_1 (x_1) = \B{h}_1\T \B{\xi}_1 = 1, \quad \dddot{\beta}_1 (x_1) = \dddot{\B{h}}_1\T \B{\xi}_1 = 0, \quad \dot{\beta}_1 (x_2) = \dot{\B{h}}_2\T \B{\xi}_1 = 0, \quad \text{and} \quad \beta_1 (x_3) = \B{h}_3\T \B{\xi}_1 = 0,
\end{equation*}
must be satisfied. Selecting, for instance, $\B{h} (x) = \{e^x, \; \sin x, \; \ln x, \; x^{-1}\}\T$, these relationships can be put in a matrix form
\begin{equation*}
    \begin{bmatrix} \B{h}_1\T \\ \dddot{\B{h}}_1\T \\ \dot{\B{h}}_2\T \\ \B{h}_3\T\end{bmatrix} \B{\xi}_1 =
    \begin{bmatrix} e^{x_1} & \sin x_1 & \ln x_1 & x_1^{-1} \\ e^{x_1} & -\cos x_1 & 2/x_1^3 & x_1^{-4} \\ e^{x_2} & \cos x_2 & 1/x_2 & -x_2^{-2} \\ e^{x_3} & \sin x_3 & \ln x_3 & x_3^{-1}\end{bmatrix} \B{\xi}_1 = \begin{Bmatrix} 1 \\ 0 \\ 0 \\ 0\end{Bmatrix}
\end{equation*}
Allowing to obtain the coefficients vector, $\B{\xi}_1$, of the $\beta_1 (x)$ function by matrix inversion. The other $\B{\xi}_i$ coefficients vectors, where $i = 2,3,4$, can be computed analogously obtaining the final following system
\begin{equation*}
    H (\B{x}) \, \Xi = \begin{bmatrix} \B{h}_1\T \\ \dddot{\B{h}}_1\T \\ \dot{\B{h}}_2\T \\ \B{h}_3\T\end{bmatrix} \begin{bmatrix} \B{\xi}_1, & \B{\xi}_2, & \B{\xi}_3, & \B{\xi}_4\end{bmatrix} = \begin{bmatrix} 1 & 0 & 0 & 0 \\ 0 & 1 & 0 & 0 \\ 0& 0 & 1 & 0 \\ 0 & 0 & 0 & 1\end{bmatrix} = I_{4\times 4}
\end{equation*}
Therefore, the inversion of the coefficients matrix, $\Xi$, provides all the $\B{\xi}_k$ coefficient vectors,
\begin{equation*}
    \Xi (\B{x}) = \begin{bmatrix} \B{\xi}_1, & \B{\xi}_2, & \B{\xi}_3, & \B{\xi}_4\end{bmatrix} = H^{-1} (\B{x})
\end{equation*}
and the $\beta_k (x, \B{x})$ polynomials are
\begin{equation*}
    \begin{Bmatrix} \beta_1 (x, \B{x}), & \beta_2 (x, \B{x}), & \beta_3 (x, \B{x}), & \beta_4 (x, \B{x})\end{Bmatrix} = \B{h}\T (x) \, \Xi (\B{x}) = \B{h}\T (x) \, H^{-1} (\B{x}).
\end{equation*}

\section{Relative Constraints}

Sometime constraints are not defined in an absolute way (e.g., $\dot{y} (2) = 1$) but in a relative way, as $\ddot{y} (0) = y (1)$. Constrained expressions can also be derived for relative constraints. In general, a relative constraint can be written as,
\begin{equation*}
    y_{\B{x}_i}^{(\B{d}_i)} = y_{\B{x}_j}^{(\B{d}_j)} \qquad \text{where if} \quad \B{x}_i = \B{x}_j \quad \text{then} \quad \B{d}_i \ne \B{d}_j.
\end{equation*}

\begin{figure}[ht]
    \centering\includegraphics[width=\linewidth]{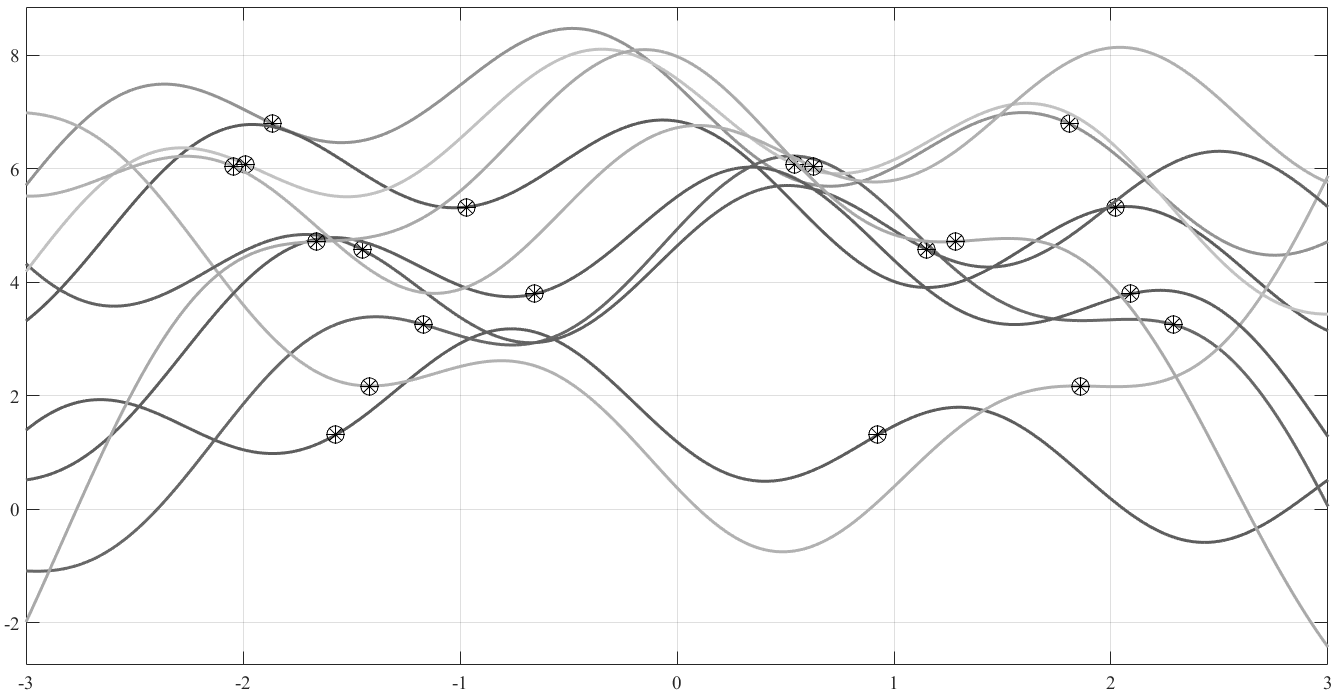}
	\caption{10 random examples satisfying $y_1 = y_2$ and $\dot{y}_1 = \dot{y}_2$}
    \label{y1uy2andyd1uyd2}
\end{figure}

To give an example, consider the problem of finding an expression satisfying the two relative constraints,
\begin{equation*}
    y_1 = y_2 \qquad \text{and} \qquad \dot{y}_1 = \dot{y}_2.
\end{equation*}
The constrained expression can be searched as,
\begin{equation*}
    y (x) = g (x) + \xi_p \, p (x) + \xi_r \, r (x),
\end{equation*}
where $p (x)$ and $r (x)$ are two assigned functions and $\xi_p$ and $\xi_r$ two unknown coefficients. The two relative constraints imply solving the system,
\begin{equation*}
    \begin{Bmatrix} g_2 - g_1\\ \dot{g}_2 - \dot{g}_1\end{Bmatrix} = \begin{bmatrix} p_1 - p_2 & r_1 - r_2 \\ \dot{p}_1 - \dot{p}_2 & \dot{r}_1 - \dot{r}_2\end{bmatrix} \begin{Bmatrix} \xi_p \\ \xi_r\end{Bmatrix},
\end{equation*}
to obtain the unknown coefficients, $\xi_p$ and $\xi_r$. This linear problem admits solution if the matrix is not singular. Matrix singularity indeed occurs if, for instance, $p (x) = 1$ and $r (x) = x$, as done in the two constraints example \#1, see Eq. (\ref{c1}). This means that the case of two relative constraints is different from the case of two absolute function constraints or one function and one derivative absolute constraint. In this relative constraint case functions $p (x)$ and $r (x)$ must admit, at least, a nonzero first derivative (in addition to span two independent function spaces). Using monomials, a potential constrained expression can be search as,
\begin{equation*}
    y (x) = g (x) + \xi_p \, x + \xi_r \, x^2.
\end{equation*}
The two constraints give the system,
\begin{equation*}
    \begin{Bmatrix} g_2 - g_1\\ \dot{g}_2 - \dot{g}_1\end{Bmatrix} = \begin{bmatrix} x_1 - x_2 & x_1^2 - x_2^2 \\ 0 & 2(x_1 - x_2)\end{bmatrix} \begin{Bmatrix} \xi_p \\ \xi_r\end{Bmatrix},
\end{equation*}
whose solution leads to,
\begin{equation*}
    y (x) = g (x) +  \dfrac{x}{x_1 - x_2} \, (g_2 - g_1) + \dfrac{x^2 - x (x_1 + x_2)}{2(x_1 - x_2)} \, (\dot{g}_2 - \dot{g}_1),
\end{equation*}
which has, again, the formal expression of Eq. (\ref{eq06}). This expression provides \emph{all} functions satisfying the two relative constraints. Figure \ref{y1uy2andyd1uyd2} provides 10 random examples using $p (x) = 1 - x^2$, $r (x) = \sin x$, and $g(x) = v + x^2 - \sin(3 x + v)$. The values of $v$, $x_1$, and $x_2$ have been randomly selected in the ranges, $v\in[0, 2\pi]$, $x_1 \in[-2.5, -0.5]$, and $x_1 \in[0.5, 2.5]$.

\subsection{Coefficients functions derivation for linear combination of relative constraints}

Consider the following constrained expression
\begin{equation}\label{GRC}
    y (x) = g (x) + \B{\xi}\T \B{h} (x).
\end{equation}
Absolute and relative constraints can be expressed by a set of $n$ linear equations
\begin{equation*}
    c_k = \B{\alpha}_k\T \, \B{y}_{\B{x}_k}^{(\B{d}_k)} = \B{\alpha}_k\T \, \B{g}_{\B{x}_k}^{(\B{d}_k)} + \B{\alpha}_k\T \, H_{\B{x}_k}^{(\B{d}_k)} \, \B{\xi}, \qquad \text{where} \qquad k = 1, \cdots, n.
\end{equation*}
where $\B{\xi}$ and $\B{\alpha}_k$ are vectors of unknown and known coefficients, respectively. For instance, for the constraint, $3 = 2 y (x_1) - \pi \ddot{y} (x_3)$, we have $c_k = 3$, $\B{\alpha}_k = \{2, \; -\pi\}$, $\B{d}_k = \{0, \; 2\}$, $\B{x}_k = \{x_1, \; x_3\}$, and $\B{y}_{\B{x}_k}^{(\B{d}_k)} = \{y_1, \; \ddot{y}_3\}$. Using a set of $n$ indefinitely differentiable basis functions (size of vectors $\B{\xi}$ and $\B{h}$), matrix $H_{\B{x}_k}^{(\B{d}_k)}$ has the expression
\begin{equation*}
    H_{\B{x}_k}^{(\B{d}_k)} = \begin{bmatrix} \B{h}_{\B{x}_{k1}}^{(\B{d}_{k_1})}, & \cdots, & \B{h}_{\B{x}_{k_m}}^{(\B{d}_{k_m})}\end{bmatrix}\T.
\end{equation*}
The $\B{\xi}$ vector can be then computed from the constraints' equations
\begin{equation*}
    \begin{Bmatrix} c_1 \\ \vdots \\ c_n\end{Bmatrix} = \begin{Bmatrix} \B{\alpha}_1\T \B{g}_{\B{x}_1}^{(\B{d}_1)} \\ \vdots \\ \B{\alpha}_n\T \B{g}_{\B{x}_n}^{(\B{d}_n)}\end{Bmatrix} + \begin{bmatrix} \B{\alpha}_1\T \, H_{\B{x}_1}^{(\B{d}_1)} \\ \vdots \\ \B{\alpha}_n\T \, H_{\B{x}_n}^{(\B{d}_n)}\end{bmatrix} \begin{Bmatrix} \xi_1 \\ \vdots \\ \xi_n\end{Bmatrix}
\end{equation*}
from which the solution is
\begin{equation*}
    \begin{Bmatrix} \xi_1 \\ \vdots \\ \xi_n\end{Bmatrix} = \begin{bmatrix} \B{\alpha}_1\T \, H_{\B{x}_1}^{(\B{d}_1)} \\ \vdots \\ \B{\alpha}_n\T \, H_{\B{x}_n}^{(\B{d}_n)}\end{bmatrix}^{-1} \begin{Bmatrix} c_1 - \B{\alpha}_1\T \B{g}_{\B{x}_1}^{(\B{d}_1)} \\ \vdots \\ c_n - \B{\alpha}_n\T \B{g}_{\B{x}_n}^{(\B{d}_n)}\end{Bmatrix}
\end{equation*}
substituting in Eq. (\ref{GRC})
\begin{equation*}
    y (x) = g (x) + \B{h}\T (x) \begin{bmatrix} \B{\alpha}_1\T \, H_{\B{x}_1}^{(\B{d}_1)} \\ \vdots \\ \B{\alpha}_n\T \, H_{\B{x}_n}^{(\B{d}_n)}\end{bmatrix}^{-1} \begin{Bmatrix} c_1 - \B{\alpha}_1\T \B{g}_{\B{x}_1}^{(\B{d}_1)} \\ \vdots \\ c_n - \B{\alpha}_n\T \B{g}_{\B{x}_n}^{(\B{d}_n)}\end{Bmatrix} = g (x) + \B{\beta}\T (x, \B{x}) \begin{Bmatrix} c_1 - \B{\alpha}_1\T \B{g}_{\B{x}_1}^{(\B{d}_1)} \\ \vdots \\ c_n - \B{\alpha}_n\T \B{g}_{\B{x}_n}^{(\B{d}_n)}\end{Bmatrix}
\end{equation*}
the expressions for the $\beta_k (x, \B{x})$ functions, is obtained
\begin{equation}
    \B{\beta}\T (x, \B{x}) = \B{h}\T (x) \begin{bmatrix} \B{\alpha}_1\T \, H_{\B{x}_1}^{(\B{d}_1)} \\ \vdots \\ \B{\alpha}_n\T \, H_{\B{x}_n}^{(\B{d}_n)}\end{bmatrix}^{-1}
\end{equation}
Therefore, the case of linear relative constraint can be expressed as
\begin{equation}\label{eq99}
    y (x) = g (x) + \ds\sum_{k = 1}^n \beta_k (x, \B{x}) \left[c_k - \B{\alpha}_k\T \, \B{g}_{\B{x}_k}^{(\B{d}_k)}\right] \quad \text{where} \quad \beta_k^{(\B{d}_k)} (x_i, \B{x}) = \delta_{ki},
\end{equation}
This equation generalizes Eq. (\ref{eq06}), where the $\beta_k (x, \B{x})$ functions are multiplying the relative constraints written in term of the function $g (x)$.

\subsection{Example of two linear combination of relative constraints}

Let's give a numerical example. Consider to find all functions satisfying the following two relative constraints
\begin{equation}\label{cnst}
    3 = 2 y (x_1) - \pi \ddot{y} (x_3) \qquad \text{and} \qquad \pi = e \dot{y} (x_1) + y (x_2) - 3 \dot{y} (x_3)
\end{equation}
where $x_1 =-1$, $x_2 = 1$, and $x_3 = 2$. Therefore, we have,
\begin{equation}
    \left\{\begin{array}{l} c_1 = 3 \\ \B{\alpha}_1 = \{2, \; -\pi\}\T \\ \B{d}_1 = \{0, \; 2\}\T \\ \B{x}_1 = \{-1, \; 2\}\T\end{array}\right. \qquad \text{and} \qquad \left\{\begin{array}{l} c_2 = \pi \\ \B{\alpha}_2 = \{e, \; 1, \; -3\}\T \\ \B{d}_2 = \{1, \; 0, \; 1\}\T \\ \B{x}_2 = \{-1, \; 1, \; 2\}\T\end{array}\right.
\end{equation}
Consider the coefficients functions selection,
\begin{equation}
    \B{h} = \{e^x, \; \sin x\}\T, \quad \to \quad \dot{\B{h}} = \{e^x, \; \cos x\}\T, \quad \to \quad \ddot{\B{h}} = \{e^x, \; -\sin x\}\T
\end{equation}
Then the $\beta_i$ functions are
\begin{equation*}
    \B{\beta}\T = \B{h}\T (x) \begin{bmatrix} \B{\alpha}_1\T \, H_{\B{x}_1}^{(\B{d}_1)} \\ \B{\alpha}_2\T \, H_{\B{x}_2}^{(\B{d}_2)}\end{bmatrix}^{-1}
\end{equation*}
where
\begin{equation*}
     H_{\B{x}_1}^{(\B{d}_1)} = \begin{bmatrix} e^{-1} & \sin(-1) \\ e^{2} & -\sin(2)\end{bmatrix} \qquad \text{and} \qquad H_{\B{x}_2}^{(\B{d}_2)} = \begin{bmatrix} e^{-1} & \cos(-1) \\ e^{1} & \sin(1) \\ e^{2} & \cos(2)\end{bmatrix}
\end{equation*}
that is,
\begin{equation*}
    \B{\beta}\T = \{e^x, \; \sin x\} \begin{bmatrix} \{2, \; -\pi\}\T \begin{bmatrix} e^{-1} & \sin(-1) \\ e^{2} & -\sin(2)\end{bmatrix} \\ \{e, \; 1, \; -3\}\T \begin{bmatrix} e^{-1} & \cos(-1) \\ e^{1} & \sin(1) \\ e^{2} & \cos(2)\end{bmatrix}\end{bmatrix}^{-1}
\end{equation*}
and, finally,
\begin{equation*}
    \B{\beta}\T = \{e^x, \sin x\} \begin{bmatrix} 2 e^{-1} - \pi e^{2} & 2 \sin(-1) + \pi \sin(2) \\ 1 +  e - 3 e^{2} & e \cos(-1) + \sin(1) - 3 \cos(2)\end{bmatrix}^{-1} \approx \{e^x, \sin x\} \begin{bmatrix} -0.0610 & 0.0201 \\ -0.3163 & 0.3853\end{bmatrix}
\end{equation*}
Therefore, \emph{all} explicit functions satisfying the constraints given in Eq. (\ref{cnst}) can be expressed by
\begin{equation}
    \begin{split} y (x) =& \, g (x) - (0.0610 \, e^x + 0.3163 \, \sin x)(3 - 2 g_1 + \pi \ddot{g}_3) + \\ ~& \qquad + (0.0201 \, e^x + 0.3853 \, \sin x) (\pi - e \dot{g}_1 - g_2 + 3 \dot{g}_3)\end{split}
\end{equation}

\section{Periodic functions}

\begin{figure}[ht]
    \centering\includegraphics[width=\linewidth]{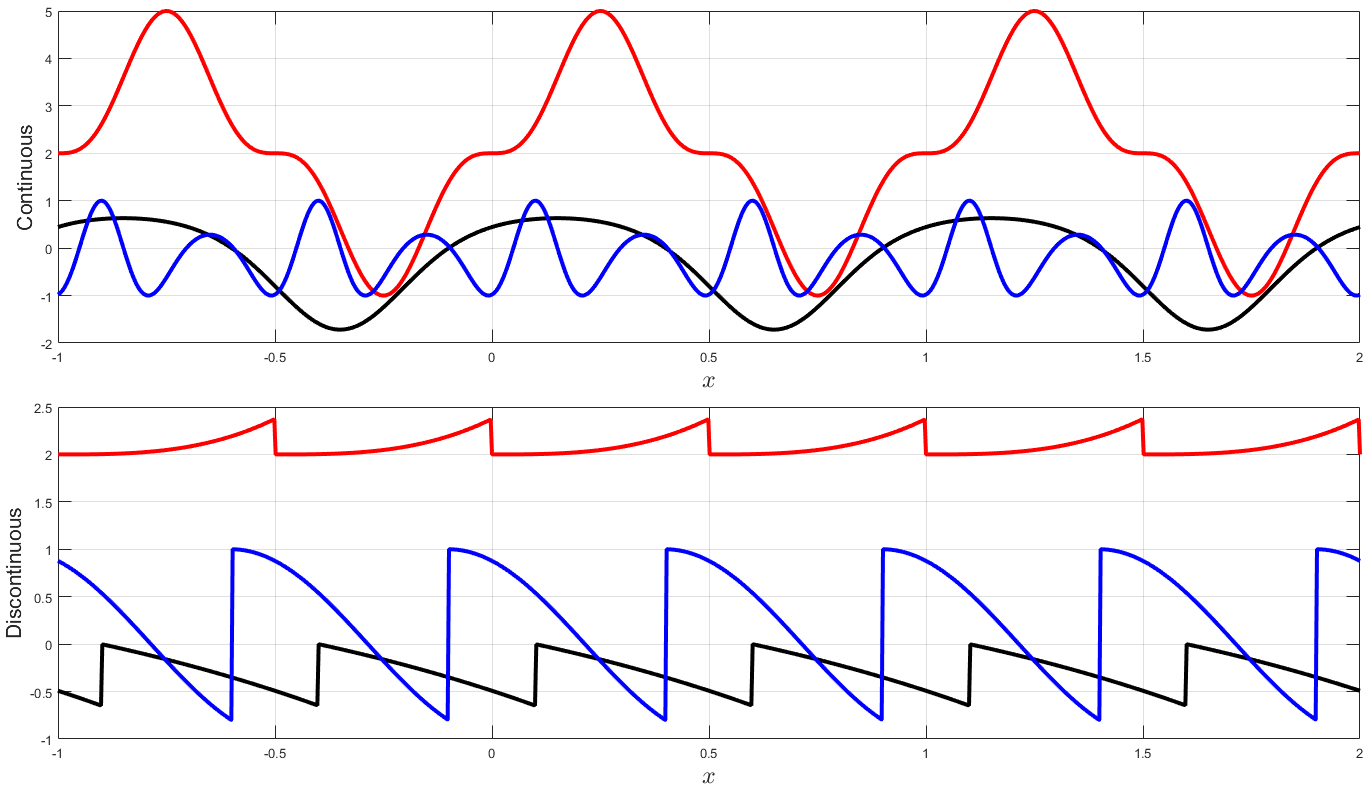}
	\caption{Example of continuous and discontinuous periodic functions}
    \label{periodic1}
\end{figure}

Periodic functions with period $T$ can be provided in continuous and discontinuous forms,
\begin{eqnarray}
    y_c (x) =& \hspace{-0.3cm} g [\Psi(x-\delta x_c, T)] \qquad & \text{(continuous)} \label{eqc} \\
    y_d (x) =& \hspace{-0.2cm} g [(x - \delta x_d) \, \text{mod} \, T]\qquad & \text{(discontinuous)} \label{eqd}
\end{eqnarray}
where $\Psi(x-\delta x_c, T)$ can be any periodic function with period $T$ (e.g., trigonometric functions), $\delta x_c$ and $\delta x_d$ can be \emph{any} constant (shifting parameters), and $g (\bullet)$ can be \emph{any} function. Figure \ref{periodic1} shows three examples of using the expressions provided by Eq. (\ref{eqc}) and Eq. (\ref{eqd}) with period $T = 0.5$, $\delta x_c = 0.4$, and $\delta x_d = 0.6$, and $\Psi(x-\delta x_c, T) \equiv \sin[\pi(x-\delta x_c)/T]$. These three functions are associated to: $g (x) = 1 - e^x$ (black), $g (x) = 2 + 3 \, x^3$ (red), and $g (x) = \cos (5 \, x)$ (blue), respectively.

The top and bottom periodic functions given in Fig. \ref{periodic1} use the continuous and the discontinuous forms, Eq. (\ref{eqc}) and Eq. (\ref{eqd}), respectively.

All periodic functions passing through the point $[x_k, y_k]$, can be expressed using the additive expression form given in Eq. (\ref{eq01}), $y (x) = g (x) + (y_k - g_k)$, that is,
\begin{equation}\label{ycdk}
    \left\{\begin{array}{ll}
    y_{ck} (x) = g [\Psi(x - \delta x_c, T)] + \{ y_k - g [\Psi(x_k - \delta x_c, T)]\} & \text{(continuous)} \\
    y_{dk} (x) = g [(x - \delta x_d) \, \text{mod} \, T] + \{y_k - g [(x_k - \delta x_d) \, \text{mod} \, T]\} & \text{(discontinuous)}
    \end{array}\right.
\end{equation}
Then, these expressions can be used along the general Waring polynomial form given in Eq. (\ref{eqw1}) to obtain \emph{all} periodic functions passing through a set of $n$ points. By doing that, the periodicity is over a line ($n = 2$), a quadratic ($n = 3$), a cubic ($n = 4$), and so on, functions.
\begin{figure}[ht]
    \centering\includegraphics[width=\linewidth]{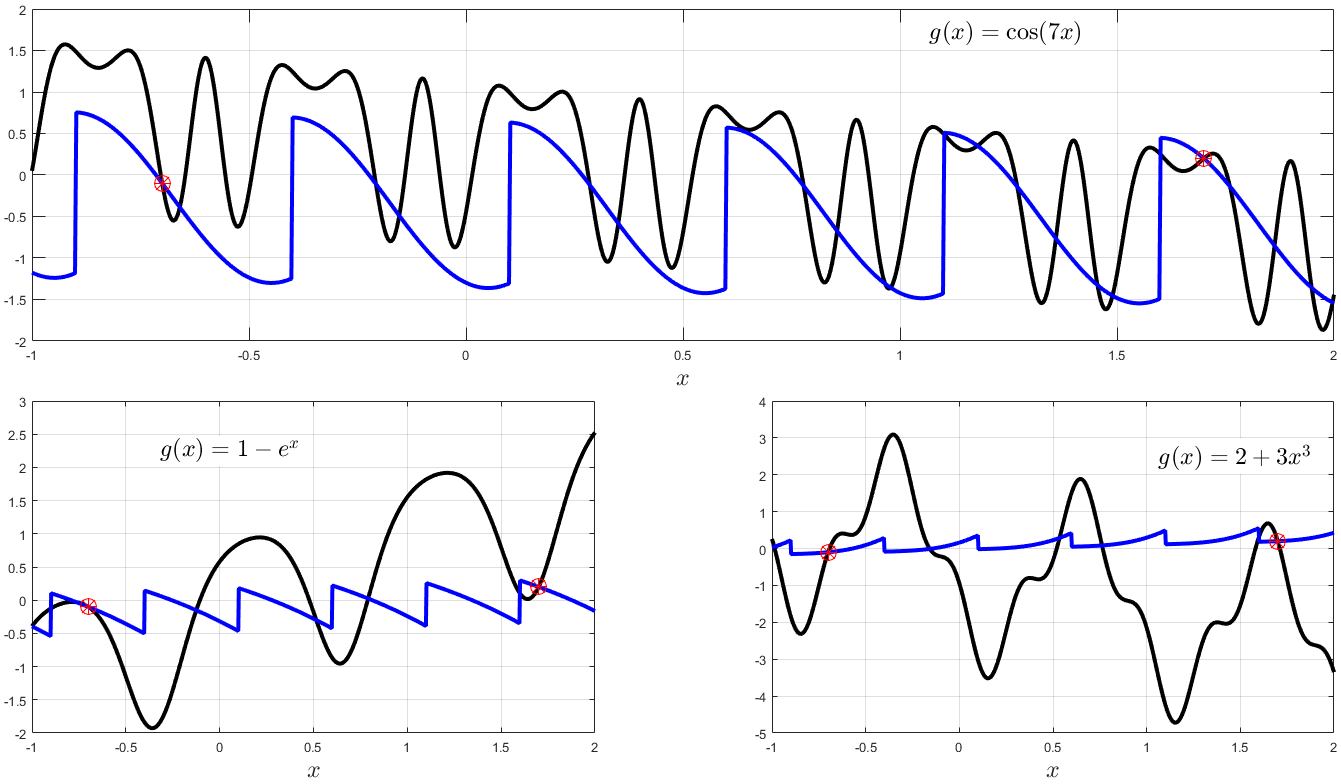}
	\caption{Example of continuous and discontinuous periodic functions}
    \label{periodic2}
\end{figure}

For instance, using two points, $[-0.7, \, -0.1]$  and $[1.7, \, 0.2]$ (indicated by red markers in all three plots of Fig. \ref{periodic2}), we can obtain the continuous (black curves),
\begin{equation*}
    y_c (x) =  \dfrac{x - x_2}{x_1 - x_2} \, y_{c1} (x) + \dfrac{x - x_1}{x_2 - x_1} \, y_{c2} (x)
\end{equation*}
and the discontinuous (blue curves),
\begin{equation*}
    y_d (x) =  \dfrac{x - x_2}{x_1 - x_2} \, y_{d1} (x) + \dfrac{x - x_1}{x_2 - x_1} \, y_{d2} (x)
\end{equation*}
for the three different expressions: $g (x) = \cos(7 \, x)$, $g (x) = 1 - e^x$, and $g (x) = 2 + 3 \, x^3$, respectively. These plots clearly show how the periodicity is over a line.

Using $n$ points the general Waring interpolation form continuous or discontinuous periodic functions over $n-1$ polynomials functions,
\begin{equation*}
    y_c (x) = \ds\sum_{k=1}^n y_{ck} \ds\prod_{i\in[1,n]}^{i\ne k} \dfrac{x - x_i}{x_k - x_i} \qquad \text{and} \qquad y_d (x) = \ds\sum_{k=1} y_{dk} \ds\prod_{i\in[1,n]}^{i\ne k} \dfrac{x - x_i}{x_k - x_i}
\end{equation*}
can be obtained, where the $n$ values of $y_{ck}$ and $y_{dk}$ are computed using Eq. (\ref{ycdk}), respectively.

\section*{Conclusions}

This study introduces \emph{constrained expressions}, providing all functions requiring to satisfy a set of assigned constraints. These expressions are not unique and are defined in terms of a new function, $g (x)$, and/or its derivatives, evaluated at constraints' coordinates. Constrained expressions can be introduced for relative constraints as well as for a set of linear relative constraints made of function and or derivatives absolute and relative constraints. Finally, the theory has been applied to periodic functions, continuous and discontinuous, with assigned period and subject to pass through a set of points.

Applications of the theory presented can be found in several fields of basic scientific research involving optimization. The fundamental property of this theory is to provide expressions that restrict the search space to just those satisfying the problem constraints. An immediate application of these constrained expressions is in solving linear nonhomogeneous differential equations with nonconstant coefficients for both, initial and boundary value problems. Reference \cite{Mortari} shows how, using constrained expressions, a least-square solution can be obtained for these two problems.

This paper introduces the first part of the \emph{theory of connections}, the one showing how to connect points and/or derivatives. A second and third part, to appear soon, will show how to connect functions and differential equations, respectively.

\end{document}